\documentclass[final]{siamltex}
\usepackage{amsmath}
  \usepackage{paralist}
  \usepackage{graphics} 
  \usepackage{epsfig} 
\usepackage{graphicx}  \usepackage{epstopdf}
 \usepackage[colorlinks = true, pdfstartview = FitV, linkcolor = blue, citecolor = blue, urlcolor = blue]{hyperref}
 \usepackage{amssymb}
\usepackage{latexsym}
\usepackage{subfigure}
\usepackage{crop}

\usepackage{color}
\definecolor{redcolor}{rgb}{1.0,0,0}
\definecolor{bluecolor}{rgb}{0,0,1.0}
\definecolor{browncolor}{rgb}{.6,.2,.2}
\newcommand{\uri}[1]{{\color{black}{}{#1}{}}}

\hypersetup{urlcolor=blue, citecolor=red}

\newcommand {\uu}  { {\bf u} }

\newcommand {\xx}  { {\bf x} }

\newcommand {\yy}  { {\bf y} }

\newcommand {\rr}  { {\bf r} }

\newcommand {\vv}  { {\bf v} }

\newcommand {\ff}  { {\bf f} }
\newcommand {\bb}  { {\bf b} }

\newcommand {\ee}  { {\bf e} }

\newcommand {\Oh} { { \mathcal O} }

\newcommand {\gb}  { {\bf g} }

\title{Discrete Processes and their Continuous Limits}

\author{Uri M. Ascher\thanks{Department of Computer Science, University of British Columbia (\tt{ascher@cs.ubc.ca}). This paper has
appeared in J. Dynamics and Games 7(2), pp 123-140, 2020}
}

\begin{document}
\maketitle




\begin{abstract}
The possibility that a discrete process can be fruitfully
approximated by a continuous one, 
with the latter involving a differential system, is fascinating. Important theoretical
insights, as well as significant computational efficiency gains may lie in store. A great success
story in this regard are the Navier-Stokes equations, which model many phenomena in 
fluid flow rather well. Recent years saw many attempts to formulate more such continuous limits,
and thus harvest theoretical and practical advantages, in diverse areas including mathematical biology, 
\uri{economics}, computational optimization,
image processing, game theory,  and machine learning.

Caution must be applied as well, however. In fact, it is often the case that the given
discrete process is richer in possibilities than its continuous differential system limit, and
that a further study of the discrete process is practically rewarding. 
Furthermore, there are situations where the continuous limit process may
provide important qualitative, but not quantitative, information about the actual discrete
process. This paper considers several case studies of such continuous limits
and  demonstrates success as well as cause for caution. 
Consequences are discussed.

\end{abstract}

\begin{keywords}numerical methods, inverse problems, differential equations, optimization, regularization\end{keywords}

\begin{AMS}{Primary: 65K05, 65L04; Secondary: 68U10}\end{AMS}

\section{Introduction}
\label{sec:intro}

The quest for continuous models that would govern in some fruitful sense families of discrete processes  is 
age-old and continuing to be fascinating~\cite{wanner10}. 
Indeed, there is much to be gained if a complicated discrete process can
be adequately approximated by a continuous one in some limit sense: 
\begin{itemize}
\item
A simpler structure more amenable to analysis may be considered, 
adding to general phenomenological understanding.
\item
Efficient numerical approaches may be constructed for the discrete problem
by going through the continuous one.
\item
Perhaps the more complicated and less elegant discrete structure may not 
even have to be considered.
\end{itemize}
A great success story in this regard are the Navier-Stokes equations: for many 
fluid flow phenomena \uri{these partial differential equations}
capture the observed flow well enough that we can forget about their derivation while concentrating on understanding their
solutions and on approximating them numerically~\cite{chorin79,chorin}.
At the same time we note, ironically, that the discrete, molecular 
description of fluid flow involves a huge number of particles, and that not all fluid phenomena are covered by this continuous model.  

A lot of exciting  
work that is relevant in our context has been carried out in recent years.
Here is but a partial list of areas and instances of interest:
 \begin{itemize}
\item
{\em Artificial time regularization}.
We return to this below, in Sections~\ref{sec:simple}, \ref{sec:reg}, and~\ref{sec:lsd}; 
see~\cite{ashudo,zangwill69} for a wider survey.
\item
{\em One-step iterative methods posed as a first order ODE discretization}.
This is a special case of artificial time continuation. Indeed, any iterative algorithm of the form
\[ \xx_{k+1} = \xx_k + \alpha_k \uri{\gb} (\xx_k) \]
which updates a current iterate $\xx_k$ (for a nonnegative iteration counter $k$) using a step size $\alpha_k > 0$
to obtain the next approximation $\xx_{k+1}$
can be viewed as a forward Euler discretization of the ordinary differential equation (ODE)
\[ \frac{d\xx}{dt} = \uri{\gb} (\xx (t)) .\]
This differential equation is obtained from the given discrete process by letting $\alpha_k \rightarrow 0$,
though, while in practice we may well want to keep the step size quite far from zero.
Herein lies a difference between qualitative and quantitative information.    
\item
{\em Two-step iterative methods posed as a second order ODE discretization}.
There has been a lot of recent activity in this regard in attempts to better understand accelerated
gradient descent methods for unconstrained optimization~\cite{polyak64,beck2017,su2014,su2016,wibisono2016,betancourt2018,zhang18,dior19}.
We will not dwell further in this article on this mushrooming area of research,
beyond what is exposed in Section~\ref{sec:lsd}.
\uri{We do believe nonetheless that the present paper is relevant in this context.} 
\item
{\em Regularizing image processing problems using 
\uri{a penalty on the discretized  gradient of an image $u$ in its continuous limit $\int_\Omega  | \nabla u|^s$},
with $s = 2$ or $s=1$.} See, for instance,~\cite{chsh,chpo2016,ashu18}.
We return to this below, in Section~\ref{sec:paradigm}.
\item
{\em Continuation methods for nonlinear equations; homotopy path, etc.} 
See, for instance,~\cite{deuflhard04,ashudo} and references therein. 
Such methods can be viewed  as artificial time integration, but again we focus on other directions here. 
\item
{\em Semiconductor equations}.
See, for instance, the text~\cite{mrs90}. A significant volume of work preceded and followed this book.
In our present context we observe that, unlike the case for fluid flow, more than one continuous process has
been derived and the practical importance of the underlying particle process has not diminished over the years.  
\item
{\em Hamilton-Jacobi and mean field in game theory}.
See, for instance,~\cite{ll07,hmc06} and many following articles such as~\cite{gomez13,burger13,gomez14} and others.
Here we do not focus on this large and relatively recent volume of work.

\item
{\em Deep learning}.
The tremendous potential in machine learning and in particular deep learning (DL) techniques has riveted scientists
and engineers in recent years~\cite{Goodfellow-et-al-2016}.
Noting that the connection between consecutive layers in a DL network  resembles a finite difference scheme,
models depending on corresponding limit differential equations and algebraic multigrid methods have been proposed
in order to better design such neural networks.
See, for instance,~\cite{ruha19}.
Interest in two-step optimization algorithms mentioned above is also related to optimization methods in neural networks.

\end{itemize}

With all the promise and excitement of the continuous-limit models, it is also important
to emphasize that care should be taken to ensure that, {\em for the given task},  all important and interesting properties of
the discrete process are captured by the continuous one.
Otherwise, there is the possibility of being restricted to a non-optimal path. Furthermore, 
it is important, and occasionally crucial, to distinguish between 
{\em qualitative} and {\em quantitative}  relationships between the discrete and the continuous.
Undoubtedly, many researchers have arrived at such a practical conclusion in specific circumstances \uri{(for instance in {\em collective dynamics})},
\uri{with 
graduate students who are responsible for the implementation of the relevant algorithms} 
often being the first to realize this.

Below we describe several short case studies that aim to highlight various aspects of the issues involved.
We start in Section~\ref{sec:simple} with a simple example which demonstrates that viewing a discrete algorithm as
a discretization of a continuous process may simply be a matter of taste, depending on what one is used to
and feels comfortable with.
Then in Section~\ref{sec:reg} we quickly derive a useful and well-known class of regularization methods for inverse problems
following the artificial time reasoning, so here is a case where the continuous limit offers
a useful and somewhat different perspective.
In Section~\ref{sec:mid}  we then turn the heat up a notch and discuss the potential, if rare, instability of
the classical implicit midpoint method for stiff ODEs. Here we see an instance where, on one hand,
the meaningful continuous limit is not what comes to mind without thought, and on the other hand,
the correct limit leads to a rigorous characterization of unstable scenarios.
In the longer Section~\ref{sec:lsd} we develop a scenario where a discretization of the simplest heat equation with constant
time step \uri{size} exhibits the usual stability limit on that step, arising from a semi-discretization in space;
but  with the freedom of choosing
a variable time step \uri{size} such a stability connection to spatial discretization does not necessarily hold.
We then connect this scenario to the gradient descent method of Section~\ref{sec:simple} and demonstrate
that the lagged steepest descent method~\cite{babo,doeasc} produces such large step sizes.
In Section~\ref{sec:paradigm} we give a quick overview of our past efforts to apply a Tikhonov-type regularization
involving differential equations, as mentioned above, to image and surface processing applications.
These efforts, and many others, have only been partially successful, and we explain why.
Section~\ref{sec:conclusions} seals the paper by offering a few general conclusions.

\uri{In each of our case studies the discrete process has a {\em step size} parameter $h$, and we examine
a limit process for it. If the step size depends on the stage or iteration $k$, denoted $h_k$, then $h = \max_k h_k$.  }    


\uri{\paragraph{\bf Remark}

The author's curiosity about the subject of the present work was first aroused before 1980.
We were constructing the general-purpose software package COLSYS~\cite{colsys81}
for boundary value ODE systems~\cite{amr}. 
Such discretized nonlinear problems often lead to systems of nonlinear algebraic equations that are difficult to solve,
and we used the damped Newton method with an automatic choice for the damping, or step size, parameter.
The ideal choice for this parameter is $=1$, where a quadratic convergence rate is achieved
for Newton's method.    
However, often in applications the Newton direction was  poor, and a rather small step size was required to obtain
decrease is some monitor objective function.
The essential difficulty was that upon moving the iterate by such a small step size the next Newton direction
was also poor, causing a repeated similar difficulty! We ended up providing means in COLSYS to switch to other
methods in such situations. But the lesson relevant  here was that, in practice, letting the step size go to $0$ was not desirable
at all, as this limit process is not ``practically smooth''. 
} 


\section{Gradient descent stability bound}
\label{sec:simple}

Consider the problem of finding the minimum of a convex $C^1$ function of $n$ variables $f (\xx)$, and denote
\begin{eqnarray}
\xx^* = \arg\min_{\xx \in \Re^n} f(\xx ) . \label{1.1}
\end{eqnarray}
The gradient descent iterative method evaluates at the current iterate $\xx_k$ the gradient $\grad f (\xx_k)$
and sets the next iterate as
\begin{eqnarray}
\xx_{k+1} = \xx_k - \alpha_k \grad f (\xx_k), \label{1.2}
\end{eqnarray}
where $\alpha_k$ is the (positive) step size and $k$ is the iteration counter, $k=0, 1, \ldots $ See, e.g.,~\cite{nw}.
   
Let us next restrict attention to the quadratic case
\begin{eqnarray}
f(\xx ) = \frac 12 \xx^TA\xx - \bb^T\xx , \label{1.3}
\end{eqnarray}
where $A$ is a given $n \times n$ real symmetric positive definite \uri{(SPD)} matrix and $\bb$ is a given real inhomogeneity vector. 
Here $\grad f(\xx) = A\xx - \bb$, and defining (for historical reasons) the residual \uri{$\rr = -\grad f(\xx) = A(\xx^* - \xx)$} we get the iteration
\begin{eqnarray}
\xx_{k+1} = \xx_k + \alpha_k \rr_k. \label{1.4}
\end{eqnarray} 
Next we fix the step size $\alpha_k = \alpha$ (called ``learning rate'' in DL parlance) and ask, 
what is the upper stability limit on $\alpha$?

There are two methods to answer this simple question:
\begin{enumerate}
\item
Observe that
for the error $\ee_k = \xx_k - \xx^*$ we have the recursion $$\ee_{k+1} = (I - \alpha A)\ee_k ,$$ 
so $\alpha$ is restricted by  $\| I-\alpha A \|_2 \leq 1$.
Since $A$ can be diagonalized by an orthogonal similarity transformation,
this bound translates to  
$1-\alpha \max \lambda_i \geq -1 ,$ where $\lambda_i > 0$ are the eigenvalues of $A$ arranged in decreasing order.
This leads to the stability bound
\begin{eqnarray}
\alpha \leq 2/\lambda_1 .   \label{1.5}
\end{eqnarray}
\item
Another way to see this is by introducing the limit ODE $\frac{d\xx}{dt} =   \rr (\xx(t))$, where $t \geq 0$ is artificial time.
Discretizing this ODE using the forward Euler method~\cite{apbook} clearly gives the gradient descent method  \eqref{1.4},
where the time step size is $h_k = \alpha_k$.
But now, we know the {\em absolute stability} bound for forward Euler! It is the bound~\eqref{1.5}.
\end{enumerate} 


\paragraph{\bf Conclusion}
Which method (for arriving at the same result) is better? This depends on one's ``comfort zone'',
or what one is used to. 
In fact, the second method
is ``obvious'' if you are used to numerical ODEs; but it is the unnecessarily longer route otherwise. 


\section{Regularizing ill-posed problems}
\label{sec:reg}

In our second case study, let us consider for simplicity the linear problem
\begin{eqnarray}
A\xx = \bb , \label{2.1}
\end{eqnarray}
where $A$ and $\bb$ are as described in Section~\ref{sec:simple}, but now $A$ is also large, sparse and ill-conditioned.
\uri{In terms of the SVD~\cite{agbook}
\begin{subequations}
\begin{eqnarray}
A = U \Sigma U^T , \label{new1a}
\end{eqnarray}
where $U$ is an orthogonal matrix and $\Sigma$ is diagonal, containing the singular values $s_i$ in decreasing
order on its main diagonal,
the solution is strictly given by $\xx = U\yy$, where the components of the vector $\yy$ are
\begin{eqnarray}
y_i = s_i^{-1} (U^T\bb)_i  \quad i = 1, 2, \ldots , n.    \label{new1b}
\end{eqnarray}
But if $s_{i_0} \ge 0$ is extremely small (representing a perturbation of $0$ due to noise, say) for some $i_0 \leq n$, then we wish to avoid \eqref{new1b}
for all $i \geq i_0$.  
The solution process then has to be regularized somehow.
\label{new1}
\end{subequations}}

The {\em truncated SVD} is one possible method, \uri{where we set $y_i = 0$ for all $i \geq i_0$.} 
But the SVD transformation involves large, full matrices,
and we want to take advantage of the sparsity of $A$. 
An appropriate Tikhonov regularization is possible, see, e.g., \cite{ehn1}; \uri{we will not describe this further here.}

An alternative to Tikhonov's method is to consider the artificial time formulation from Section~\ref{sec:simple}, \uri{written as}
\begin{eqnarray}
\frac  {d\xx}{dt} = \bb - A\xx . \label{2.2}
\end{eqnarray}
This well-known method has been called {\em exponential filtering} in~\cite{carezh,care}.
For the artificial time ODE, at time {$t$} the effect of the reciprocal singular value  {$s_i^{-1}$} is replaced by
{$\omega(s_i^2) s_i^{-1}$}, where
\begin{eqnarray}
\omega (s) = 1 - e^{-ts} . \label{2.2} 
\end{eqnarray}

Thus, the effect of small singular values gets dampened while
large ones remain almost intact for an appropriate  {\em finite time} instance $t$.
This artificial finite time serves as the regularization parameter~\cite{ashudo},
\uri{gauging the amount by which inverse singular values are dampened}.
Next, we can use {forward Euler}, or \uri{some} other time discretization hopefully with \uri{large step sizes $h_k$}, to roughly integrate the ODE.

\begin{figure}[htb]
  \begin{center}
    \includegraphics[scale=0.45]{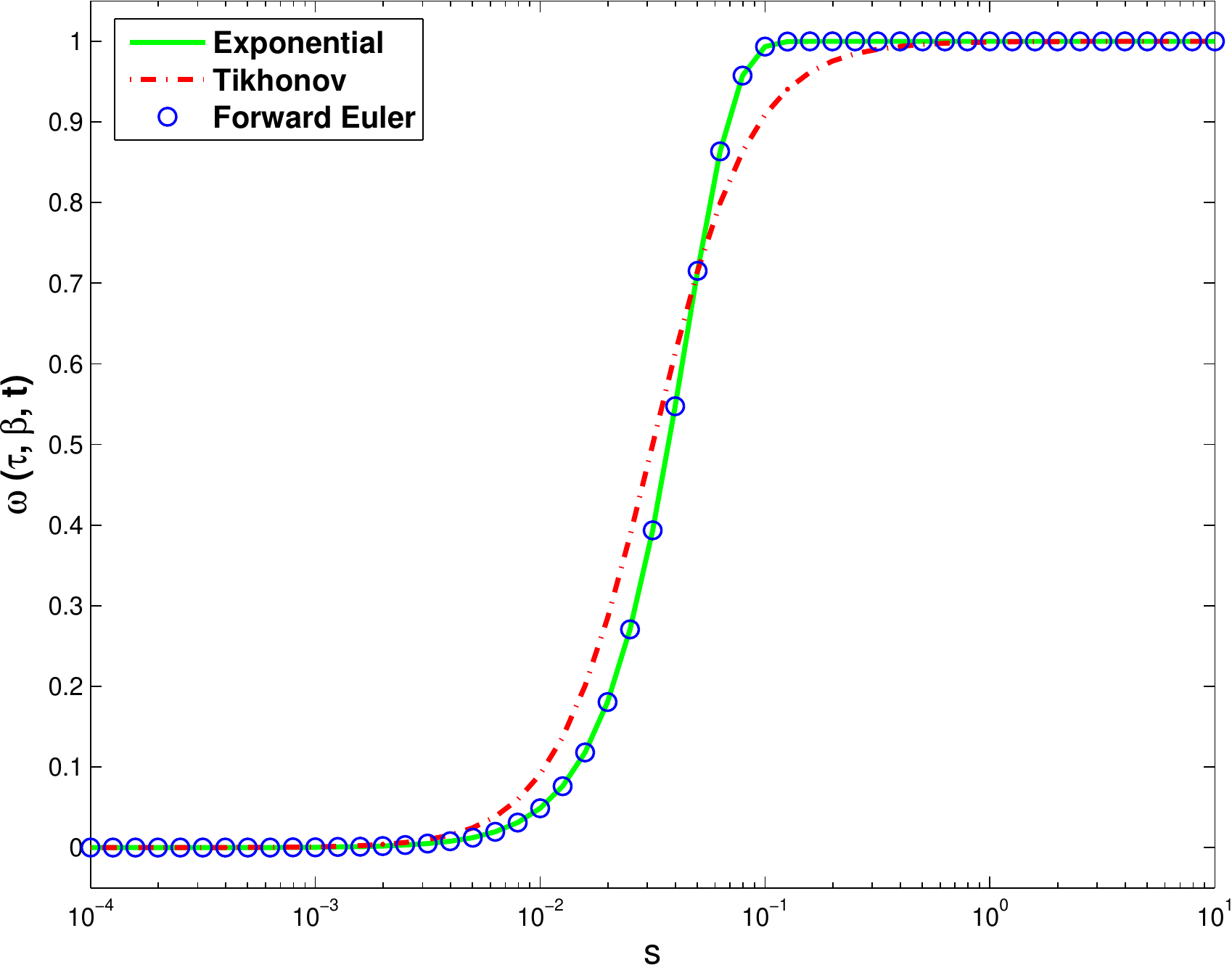}
    \end{center}
\caption{Exponential filter  $\omega (s) = 1 - e^{-ts}$
and Tikhonov filter $\omega_{\rm T}  (s) = \frac {s}{s + \beta} $
for $t\beta = 1/2$. \label{fig1}}
\end{figure}

The resulting filter function, as it turns out, behaves quite well and not so differently from the Tikhonov filter
with regularization weight $\beta$ given by
\begin{eqnarray}
\omega_{\rm T}  (s) = \frac{s}{s+\beta} . \label{2.3} 
\end{eqnarray}
For the choice of artificial time $t_\beta=\frac 1{2\beta}$ the curves in Figure~\ref{fig1} are quite similar.

Some relevant questions, discussed in~\cite{ashudo}, are:
\begin{enumerate}
\item
Using an explicit discretization method, can we take large steps even if the ODE is stiff?
\item
Is the resulting regularization effective?
(Note that we can view it  as a {\em continuation method}: we want to get to {$t=t_\beta$} quickly, without approximating the ODE trajectory accurately.)
\end{enumerate}


\paragraph{\bf Conclusion}
We will not dwell on this further, \uri{since this case study is documented in a similar context elsewhere.}
But \uri{let us}  observe that here the artificial ODE yields a practically useful regularizer!
At the same time, the task of approximating the ODE~\eqref{2.2} should not be taken
too zealously: in fact, using a higher order method such as explicit two-stage Runge-Kutta~\cite{apbook}
does not prove more useful than forward Euler for the given purpose.
\uri{Moreover, this entire regularization approach does have competition: for instance, our own favourite 
has often
been a subspace regularization method, applying just a few conjugate gradient iterations to a singular linearized system~\cite{doas1,rodoas}.}


\section{Stability of the midpoint method}
\label{sec:mid}

In our third case study we start with a given ODE system (so the independent variable $t$ is not artificial)
and consider another, implicit discretization.
Let us denote $\dot u = \frac{du}{dt}$.   
For the ODE system
\begin{eqnarray}
\dot \uu = \ff (t,\uu ), \quad 0 \leq t \leq T, \label{1}
\end{eqnarray}
the well-known implicit midpoint method reads
\begin{eqnarray}
\frac {\uu_{k+1} - \uu_k }h = \ff \left(\frac{t_{k+1}+t_k}2, \frac{\uu_{k+1}+\uu_k}2 \right) . \label{2}
\end{eqnarray}
Here  $h = h_k$  is a step  size from time $t_k$ (where for an initial value problem $\uu_k$ is known) to time $t_{k+1}$ (where $\uu_{k+1}$ is  unknown).

This method is implicit, second order accurate, A-stable, AN-stable (meaning it is stable for $f = \lambda(t)u$ for any $\Re (\lambda) \leq 0$),
 algebraically stable, symplectic for Hamiltonian systems, symmetric \uri{(reversible)}, and it conserves quadratic invariants.
 Collocation at Gaussian points generalizes midpoint  to higher order methods and ODE systems~\cite{amr}.
 The method has been applied  in the context of boundary value ODEs~\cite{amr}, geometric integration~\cite{hlw}
 and stiff initial value ODEs (IVODEs).
 We continue with the latter.
 
 Despite all its good properties (and perhaps because of some of them), the implicit midpoint method is not always stable for well-posed, stiff IVODE systems!
 Our task is next to analyze this situation.
 To recall, a stiff ODE typically has more than one time scale, and we aim to employ a time step $h$ that is
 commensurate with the slow scale but not with the fast one~\cite{hw,bawi,kale}.
 In such a constellation, though, {\em the continuous limit of the discrete structure} \eqref{2} {\em \uri{is}  not the original
 problem~\eqref{1}}. This is so, because to get from \eqref{2} to \eqref{1} we have to let $h \rightarrow 0$,
 and once we do that the problem is no longer stiff in the sense that the fast scale is eventually also resolved by the small step size.
 Instead, we have to resort to numerical singular perturbation techniques. 
 
 To study this further, consider the linear non-autonomous  ODE system
 \begin{eqnarray}
 \dot  \uu =  \omega A(t) \uu, \label{3} 
 \end{eqnarray}
 where $\omega$ is a parameter that is allowed to grow large  (so we must consider $\omega \gg 1$), 
  and $\| A \| = {\Oh} (1)$ and nonsingular. It is possible to generalize \eqref{3} to the case where $A = A(t; \omega )$.
  In \eqref{3} we are essentially studying a linearized version of~\eqref{1} that governs the propagation of an error.
 Assume that the IVODE for~\eqref{3} is well-posed (or ``stable'').
 The midpoint method reads
 \begin{eqnarray}
\frac {\uu_{k+1} - \uu_k }h =  \omega A\left(\frac{t_{k+1}+t_k}2 \right) \frac{\uu_{k+1}+\uu_k}2 , \label{4}
\end{eqnarray}
but we cannot conclude much when 
\begin{eqnarray}
\gamma = \frac 14 \omega h^2 \label{5}
\end{eqnarray}
 is fixed and not small,  because we can't take the limit of $h \rightarrow 0$ without also changing $\omega$
 and thus the ODE~\eqref{3},  as explained above. 
Specifically,  what we are looking at instead is the limit process of $\left( h \rightarrow 0, \ \omega \rightarrow \infty \right)$, such that $\gamma$ stays constant
(e.g., $\gamma = 1$).

We next introduce an old but not necessarily well-known  trick~\cite{kreiss72}.
For analysis purposes, let us 
define for each integer $k$
\begin{eqnarray}
\vv_k = (-1)^k \uu_k . \label{6}
\end{eqnarray}
Obviously, \uri{the} boundedness properties of \uri{the sequences} $\{\vv_k\}$ and $\{\uu_k\}$ are the same: $|\uu_k| = |\vv_k| ~\forall k$.

Substituting \eqref{6} in \eqref{4} and cancelling out $(-1)^k$ we obtain
 \begin{eqnarray}
\frac {\vv_{k+1} + \vv_k }h = \omega A\left(\frac{t_{k+1}+t_k}2 \right) \frac{\vv_{k+1}-\vv_k}2 . \nonumber
\end{eqnarray}
Multiplying throughout by $2/h$ \uri{then yields} an approximation for $\dot \vv$ at the right hand side.
Multiplying further by $(\omega A)^{-1}$ we can rewrite this as
 \begin{eqnarray}
\frac {\vv_{k+1} - \vv_k }h = \frac 1{\gamma} A^{-1}\left(\frac{t_{k+1}+t_k}2 \right) \frac{\vv_{k+1}+\vv_k}2 . \label{7}
\end{eqnarray} 
Now, for \eqref{7} with $\gamma$ fixed we can finally take the limit $h \rightarrow 0$, obtaining the {\em ghost ODE}
\begin{eqnarray}
\dot \vv = \gamma^{-1}  A^{-1}(t)  \vv . \label{8}
\end{eqnarray}
{\em The stability of the midpoint method therefore depends on the stability of the ghost ODE~\eqref{8}, 
not on that of the given ODE \eqref{3}}.

The exercise to design a stable IVODE problem~\eqref{3} such that the IVODE for~\eqref{8}
is unstable is solved in~\cite{a2} for the limit case $1/\omega = 0$. It can be extended to a large but finite $\omega$. 
Then the computed solution using the implicit midpoint method applied to a well-posed problem can be made to blow up;
see~\cite{a2}.\footnote{\uri{Further, although it is not strictly a matter of stability, we urge the
interested reader to  see Example~2.1 in~\cite{ar2},
where the analysis trick \eqref{6} is used to show that the symplectic midpoint method could produce consistently wrong results
for highly oscillatory non-autonomous Hamiltonian systems.}}


The above analysis can be summarized \uri{as follows}:
\begin{theorem}
Assume that the given IVODE for~\eqref{3} is well-posed, and
consider the midpoint discretization~\eqref{4} with $\uu_0 = \uu (0)$ for $k = 0, 1, \ldots, N-1$ and $Nh = T$.

Then this method is stable along the ray $\left( h \rightarrow 0,~\omega \rightarrow \infty \right)$ such that $\gamma$ defined in~\eqref{5}
is held fixed, iff
the IVODE for the ghost ODE~\eqref{8} is well-posed.  
\end{theorem}

If $A$ is a constant matrix, then the stability of the midpoint method is guaranteed. But in the more general case it is not.
Fortunately, such midpoint instability is rather rare in practice, and this is important especially for boundary value ODEs,
where general-purpose codes typically implement symmetric schemes.\\

\paragraph{\bf Conclusion}
In this case study we have learned two lessons. The first is that just because a differential equation looks like
a plausible continuous limit does not mean that it is the correct continuous limit for a discrete process
in the context of a given task.
The second lesson is that with some extra care and ingenuity it may be possible to use a continuous limit to prove
interesting theoretical results that are of practical importance.


\section{Forward Euler for the heat equation and chaotic descent}
\label{sec:lsd}

In this case study we obtain some unexpected and perhaps counter-intuitive results
by considering time-stepping with highly variable step sizes.

Consider the simple heat equation on a unit square in space and time
\begin{eqnarray}
 \uri{\frac{\partial u}{\partial t} \equiv }  u_t = \Delta u + b, \quad 0 < x,y < 1, \; t \geq 0 \label{5.1}
 \end{eqnarray}
 with homogeneous boundary conditions (BC) and initial condition $u(0,x,y) = u_0 (x,y)$ \uri{that satisfies the BC}.
 Here $\Delta u = u_{xx} + u_{yy}$ is the Laplacian, and $b > 0$ is a constant for simplicity. 
 Integrating in time to steady state, we obtain
 the model Poisson problem via a {\em continuation method} in the time variable
\begin{eqnarray}
 -\Delta u = b, \quad 0 < x,y < 1
 \label{5.2}
 \end{eqnarray} 
subject to the same BC.\footnote{\uri{A time-honoured {\em continuation method} for finding an approximte solution for a tough nonlinear
elliptic PDE $\psi(u) = 0$ is to numerically integrate the parabolic PDE $u_t = \psi(u)$, starting from an initial guess $u_0$,
until a suitable error tolerance is satisfied. For our purpose here we consider the simplest linear version of this.}}

 A forward Euler discretization of~\eqref{5.1} in time gives the recursion
 \begin{eqnarray}
 u_{k+1} = u_k + h_k \Delta u_k + h_kb , \quad k = 0, 1, \ldots  . \label{11}
 \end{eqnarray}
This is clearly unstable for any series of positive time \uri{step sizes} $h_k$.
The {\em unconditionally unstable} semi-discretization \eqref{11} is then our continuous model.
 
 Next, we discretize the partial differential equation (PDE)~\eqref{5.1}  in space first, with constant spacing $\xi = 1/\sqrt{n}$ in both $x$ and $y$ directions.
 Reshaping the 2D array of unknowns at spatial mesh points into a vector $\uu$ of length $n$ we obtain an ODE system
 \begin{eqnarray}
 \frac{d\uu}{dt} = -A \uu +\bb, \label{11.5}
 \end{eqnarray}
where the $n \times n$ matrix $A$ is SPD. Using  straightforward centred differences for the Laplacian, 
$A$ has $5$ \uri{nonconsecutive} nonzero diagonals
(see, e.g., Section~7.1 of~\cite{agbook}).

Our discrete model is obtained upon applying the same forward Euler method as for the continuous model
to \eqref{11.5}, obtaining
 \begin{eqnarray}
 \uu_{k+1} = \uu_k + h_k (\bb - A\uu_k)  . \label{12}
 \end{eqnarray}
This recusion is the same as \eqref{1.4} in Section~\ref{sec:simple}.
Evaluating eigenvalues, the stability \uri{restriction} \eqref{1.5} translates here to the condition
\begin{eqnarray}
h \leq \frac 14 \xi^2 , \quad {\rm where~} h = \max_k h_k . \label{5.5}
\end{eqnarray} 
So, when we let $\xi \rightarrow 0$ 
the stability region for $h$ shrinks like $\xi^2$.
The unconditional instability as above is effectively reached in the continuous limit for \eqref{11}.
For a uniform step size, then, the continuous model governs the discrete process well.

But next we ask,  is this shrinking stability limit on the maximum step size necessary even if we allow $h_k$ to vary?
The answer turns out to be negative! 
Indeed, consider a sequence of several small step sizes $h_k$ effectively reducing the amplitude of high eigenvalue modes of the residual,
followed by a large step reducing low mode amplitudes while inflaming the high modes to some extent.
The latter are reduced again using further small steps. 
The largest step size $h$ depends on the low modes
rather than the high modes, and as
such it can be independent of $\xi$.
This possibility is not available for the (semi-) continuous limit scheme, only the fully discrete one!


\uri{\paragraph{\bf Remark}
Before proceeding to see an actual method which produces large step sizes, let us note that it is $h$ and 
not $\xi$ which is in our focus of interest here. Indeed, the entire PDE setup can be replaced by a time-dependent linear ODE
system discretized as in~\eqref{12}, 
where the constant SPD matrix $A$ has eigenvalues that grow larger and larger
and the question is whether $h$ must then tend to zero for convergence to steady state.} 

\subsection{Faster gradient descent}
\label{sec:faster}

To actually obtain such a sequence of time steps as described above, we use the interpretation of
the forward Euler discretization as that of gradient descent, as in Section~\ref{sec:simple}.
Since we are switching from numerical PDEs to optimization,
let us change notation slightly to conform to different standards by letting  {$\xx \leftarrow \uu, \ \alpha_k \leftarrow h_k$}.
Then our forward Euler  for the discretized heat equation becomes
{gradient descent} for the quadratic objective function~\eqref{1.3}
\[ \xx_{k+1} = \xx_k + \alpha_k \rr_k, \quad \rr_k = \bb - A\xx_k  .\]
We next concentrate on choosing the step size $\alpha_k = h_k$.

The steepest descent (SD) choice is 
\begin{eqnarray}
\alpha_k^{SD} = \frac{\rr_k^T\rr_k}{\rr_k^TA\rr_k} \equiv \frac {(\rr_k, \rr_k)}{(\rr_k, A\rr_k)}  . \label{5.6}
\end{eqnarray}
This is the ``greedy choice'', obtained upon performing exact line search. It is known to perform well for the first
few steps but to become slower later on~\cite{akaike}. Asymptotically \uri{SD performs as well} 
as the best
constant step choice, for which the bound~\eqref{5.5} holds.

The lagged steepest descent (LSD) step size~\cite{babo} uses the same expression evaluated at the previous iterate:
 \begin{eqnarray}
\alpha_k^{LSD} = \frac{\rr_{k-1}^T\rr_{k-1}}{\rr_{k-1}^TA\rr_{k-1}}  =  \frac {(\rr_{k-1}, \rr_{k-1})}{(\rr_{k-1}, A\rr_{k-1})}  . \label{5.7}
\end{eqnarray}
This then is a two-step method, using information about both $\xx_k$ and $\xx_{k-1}$ to define $\xx_{k+1}$.\footnote{In~\cite{doeasc}
we found experimentally that the overall best  faster gradient descent methods are the two-step ones.}\\

\paragraph{\bf Example}
Returning to our heat equation with $b = 1$, we apply the scheme~\eqref{12} with $h_k = \alpha_k^{LSD}$ to advance to steady state,
stopping when $\| \rr_k \| \leq 10^{-6} \| \rr_0 \|\ .$
The maximum  step size $h$ is recorded in Table~\ref{tab1} for a decreasing sequence of spatial steps $\xi$.
\begin{table}[htb]
\begin{center}
\begin{tabular}{|l|cccc|}
\hline
$\xi$ & $2^{-5}$ & $2^{-6}$ & $2^{-7}$ & $2^{-8}$    \\
\hline
$h$ & .05 & .039  & .043 & .035 \\
\hline 
\end{tabular}
\caption{Maximum time step $h$ for the heat-to-Poisson \uri{process} 
as a function of spatial step $\xi$.
\label{tab1}}
\end{center}
\end{table}
\uri{Evidently, the maximum step size $h$ does not shrink as a function of $\xi$ in this example.}

Continuing our discussion of the LSD method, it is obvious from Table~\ref{tab1} that the maximum step size $h$ does not
satisfy the bound~\eqref{5.5}.
In other words, some step sizes $h_k=\alpha_k$ disobey the fixed-step stability limit which relates to the continuous PDE.
This causes a chaotic effect in the resulting method \uri{which nonetheless converges}~\cite{doeasc}.

There are many other faster gradient descent variants; our favourite among those is the one where
the SD step~\eqref{5.6} is simply updated only at every second iteration~\cite{rasv}.
In experiments its efficiency is comparable to that of LSD.

\begin{figure}[htb]
 \begin{center}
 \leavevmode
 \includegraphics[scale=0.17]{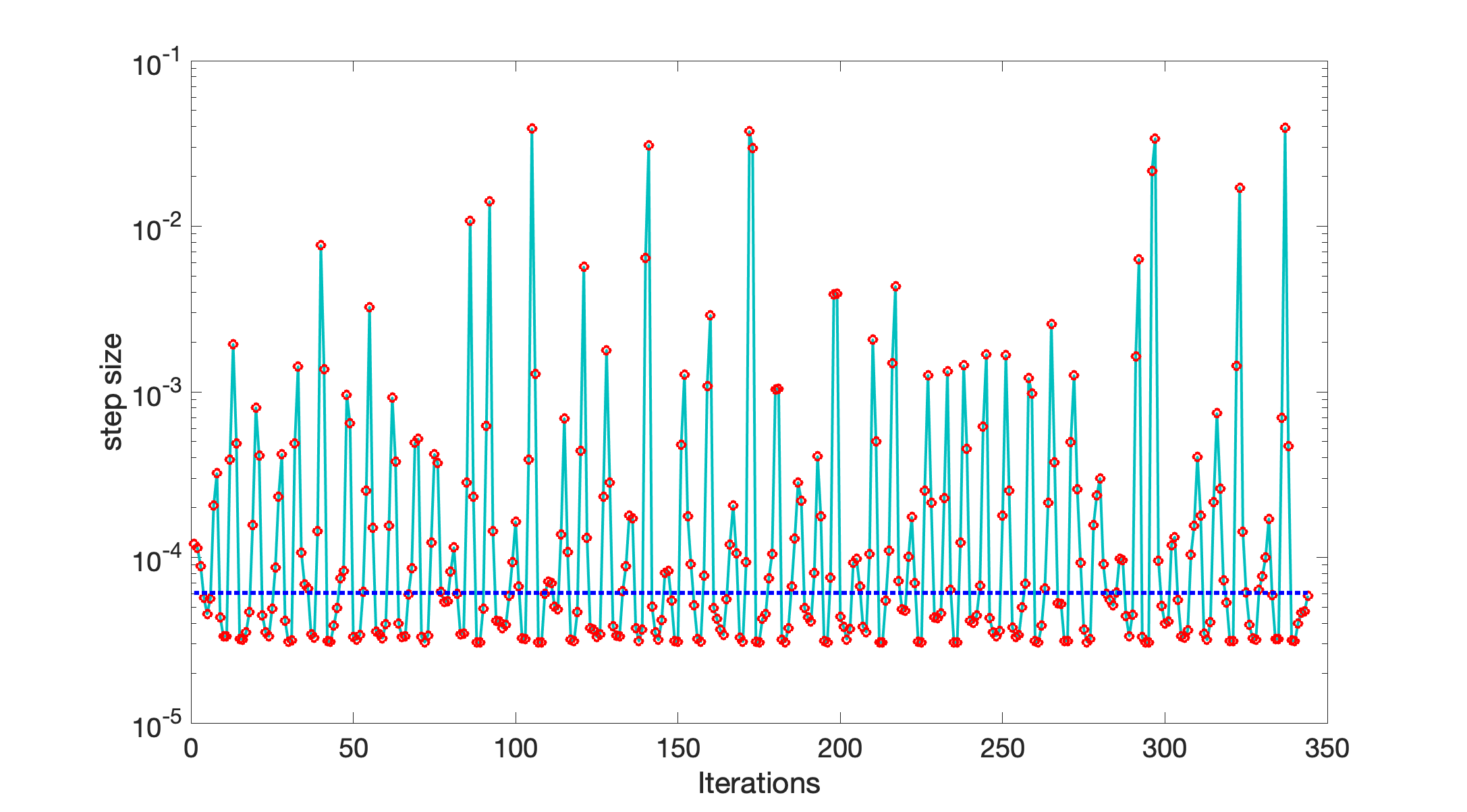}
\caption{Gradient descent with step sizes by~\eqref{5.7} for
the discretized heat equation
with $n = 63^2$. 
The calculated step sizes are displayed vs iteration counter $k$.
The stability limit for a constant step size gives the straight blue line.  \label{fig:stepsizes}}
 \end{center}
\end{figure}

Figure~\ref{fig:stepsizes} displays the step sizes observed in the LSD iteration.
Notice how much larger the largest step sizes are from the maximum constant step size.

\begin{figure}[htb]
 \begin{center}
 \leavevmode
 \includegraphics[scale=0.17]{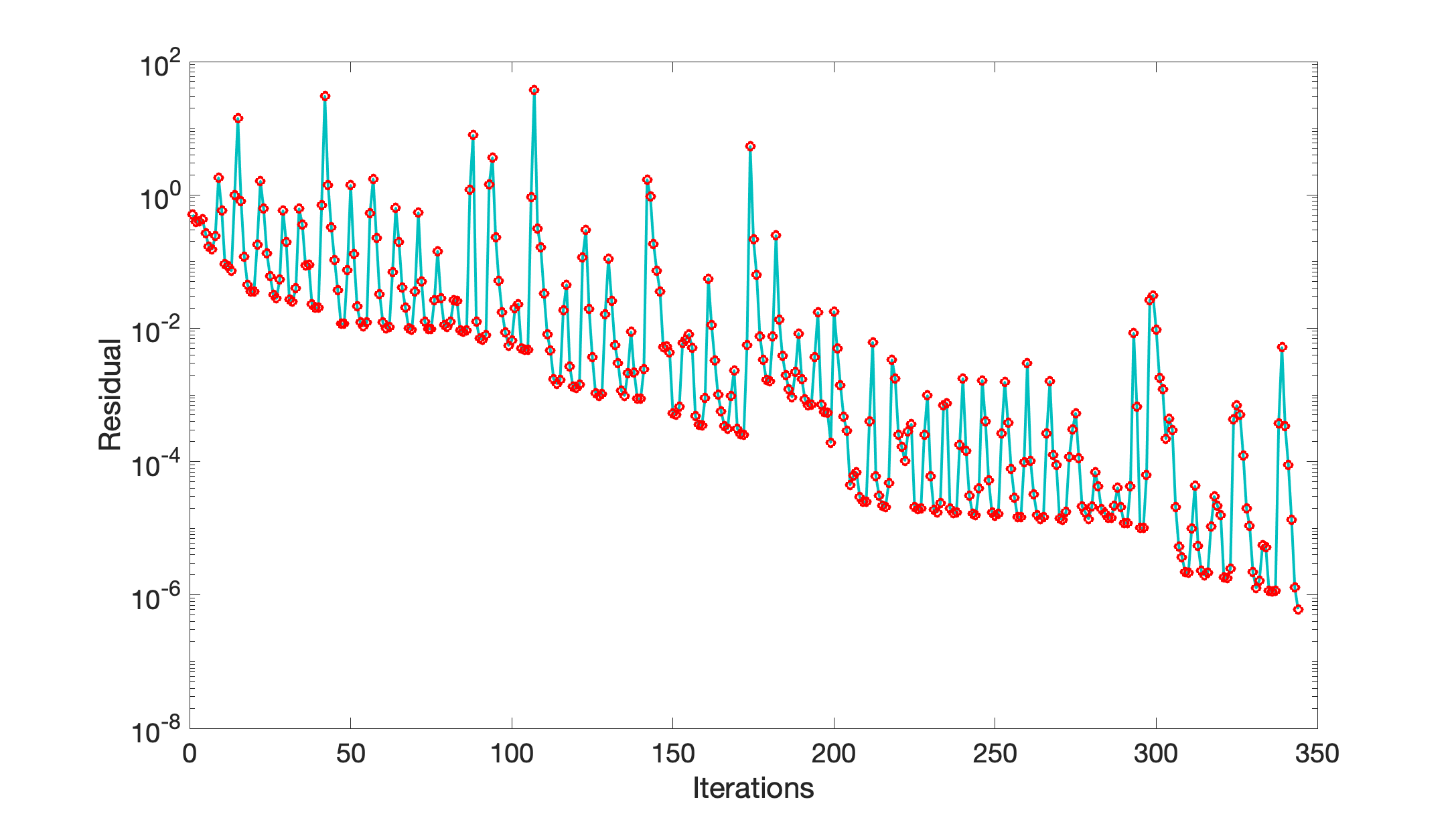}
\caption{Convergence behaviour of LSD for the model Poisson problem  with $n = 63^2$.
The errors $\| \rr_k \|$ 
are displayed as a function of iteration counter $k$. 
\label{fig:lsd_rerr}}
 \end{center}
\end{figure}

\begin{figure}[htb]
 \begin{center}
 \leavevmode
 \includegraphics[scale=0.17]{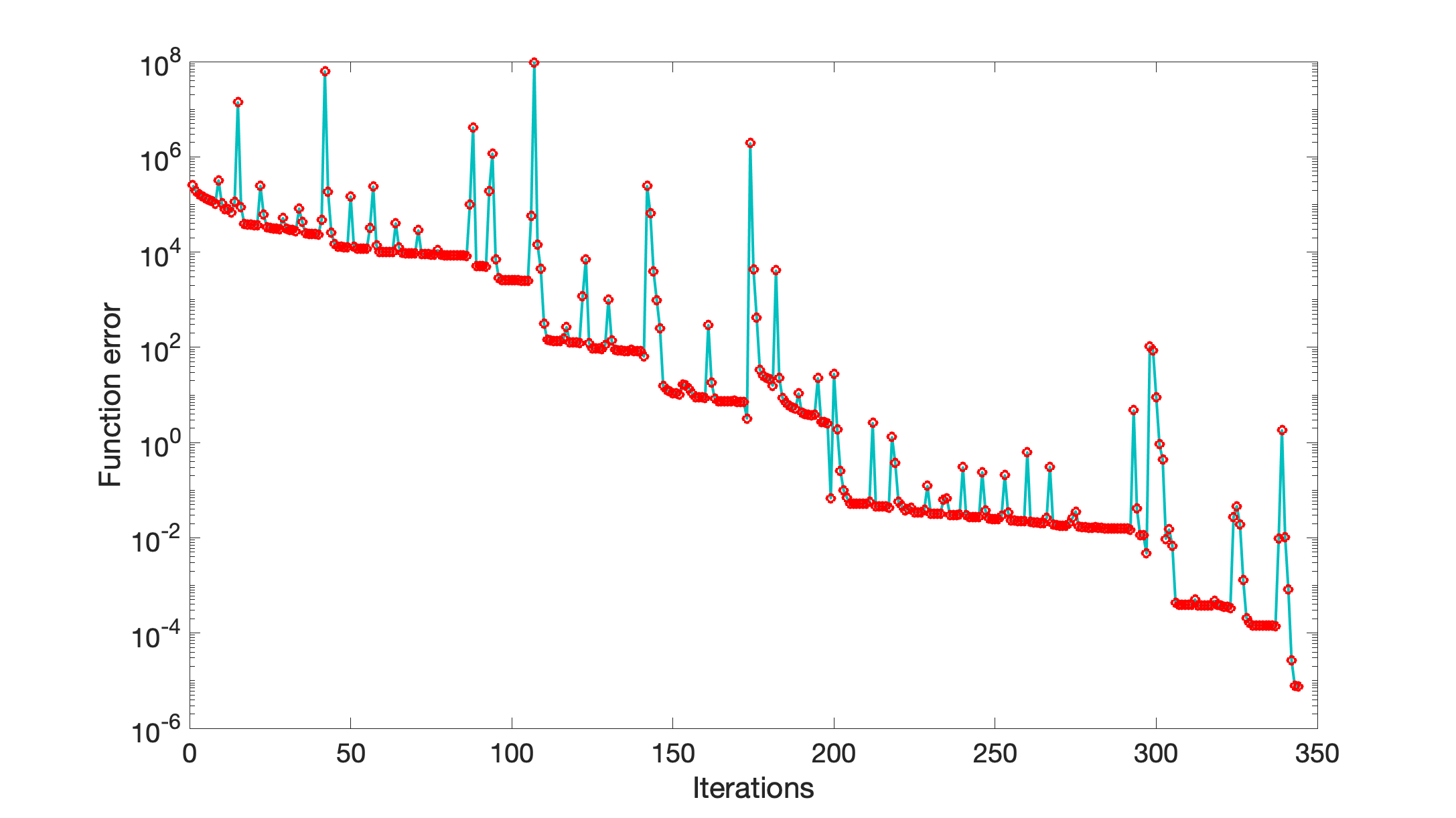}
\caption{Convergence behaviour of LSD for the model Poisson problem  with $n = 63^2$.
The errors $f( \xx_k) - f(\xx^*)$ 
are displayed as a function of iteration counter $k$. 
\label{fig:lsd_ferr}}
 \end{center}
\end{figure}

The actual residual norm $\| \rr_k \|_2$ and objective function error $f(\xx_k) - f (\xx^*)$ with $f$
given by~\eqref{1.3} are given in Figures~\ref{fig:lsd_rerr} and ~\ref{fig:lsd_ferr}, respectively.
In these figures we observe a rather non-monotone convergence, and yet the overall convergence rate 
does not appear to be slow.

As mentioned in Section~\ref{sec:intro}, there has been a lot of recent interest in  two-step optimization methods
that accelerate gradient descent.
We will not get into related ODE formulations, and instead mention the most popular of these.
The celebrated Nesterov's method~\cite{nesterov1983method} 
for the problem~\eqref{1.1}, i.e., the unconstrained minimization of a function $f(\xx)$,  is given by
\begin{subequations}
\begin{eqnarray}
\yy_{k+1} &=& \xx_k + \beta_k (\xx_k - \xx_{k-1}), \label{9a5a}\\
\xx_{k+1} &=& \yy_{k+1} - \alpha_k \grad f (\yy_{k+1})  \label{9a5b}.  
\end{eqnarray}
Note that, unlike LSD, the gradient is not evaluated at $\xx_k$, but rather at its modification $\yy_{k+1}$.
Hence, the search direction is not the residual $\rr_k$.
\label{9a5}
\end{subequations}

It can be shown that under suitable conditions this method satisfies 
\begin{eqnarray}
f(\xx_k) - f(\xx^*) = \Oh \left( k^{-2}\right) . \label{5.8}
\end{eqnarray} 
Furthermore, there is an oracle that says that such a rate is optimal.
Note that SD and constant-step gradient descent  
do not achieve \eqref{5.8}~\cite{nesterov18}.
It is therefore interesting to compare the method~\eqref{9a5} to the faster gradient descent LSD for the
case where $f$ is convex quadratic as in~\eqref{1.3}. 

However, the non-monotonic convergence displayed in Figures~\ref{fig:lsd_rerr} and~\ref{fig:lsd_ferr} makes such
a comparison difficult.
Fortunately, these errors can be ``monotonized'' by ditching the heat equation interpretation
and evaluating, following the $k$th iteration, 
\begin{subequations}
\begin{eqnarray}
\ell_k = \arg\min_{1 \leq l \leq k} \| \grad f (\xx_l)\| . \label{ella}
\end{eqnarray}
 Note that storage and CPU time can be saved by evaluating this at each step $k > 0$ as
 \begin{eqnarray}
\ell_k = \begin{cases} k & {\rm if~}  \| \grad f (\xx_k)\| \leq \| \grad f(\xx_{\ell_{k-1}}) \| \cr
             \ell_{k-1} & {\rm otherwise} \end{cases}  . \label{ellb}
\end{eqnarray}
We can then report $\xx_{\ell_k}$ at any real time, if required, as the ``best'' approximate solution at the end of iteration $k$ (although
we don't use it for the evaluation of the next iteration), because this represents  what we know at the end
of step $k$ better than $\xx_k$ does, in general. So to assess the iterative process we may record also
\begin{eqnarray}
{\rm Egrad}_k = \| \grad f(\xx_{\ell_k}) \| , \quad {\rm and~~~}   {\rm  Ef}_k =  f(\xx_{\ell_k}) - f(\xx^*) , 
\label{ellc}
\end{eqnarray}
the latter still involving a generally unknown quantity.
\label{ell}
\end{subequations}

The important point to remember is that all we need to report in practice is the computed solution upon termination. 
So, if our
termination criterion is $\| \grad f (\xx_k) \| \leq {\tt tol}$ for some error tolerance {\tt tol}, then upon
achieving this we stop. 

For the quadratic optimization problem we use the notation
\begin{eqnarray}
{\rm Er}_k = \| \rr_{\ell_k} \| .
\label{errs}
\end{eqnarray}
These residual error measures are by definition monotonically non-increasing,  and calculating them at the end of each iteration
is straightforward.
Furthermore, stopping the iteration at the first $k$ such that $\| \rr_k \| \leq {\tt tol} \| \bb \|$ is the same as
using  ${\rm Er}_k \leq {\tt tol} \| \bb \|$ for termination.

\begin{figure}[htb]
 \begin{center}
 \leavevmode
 \includegraphics[scale=0.17]{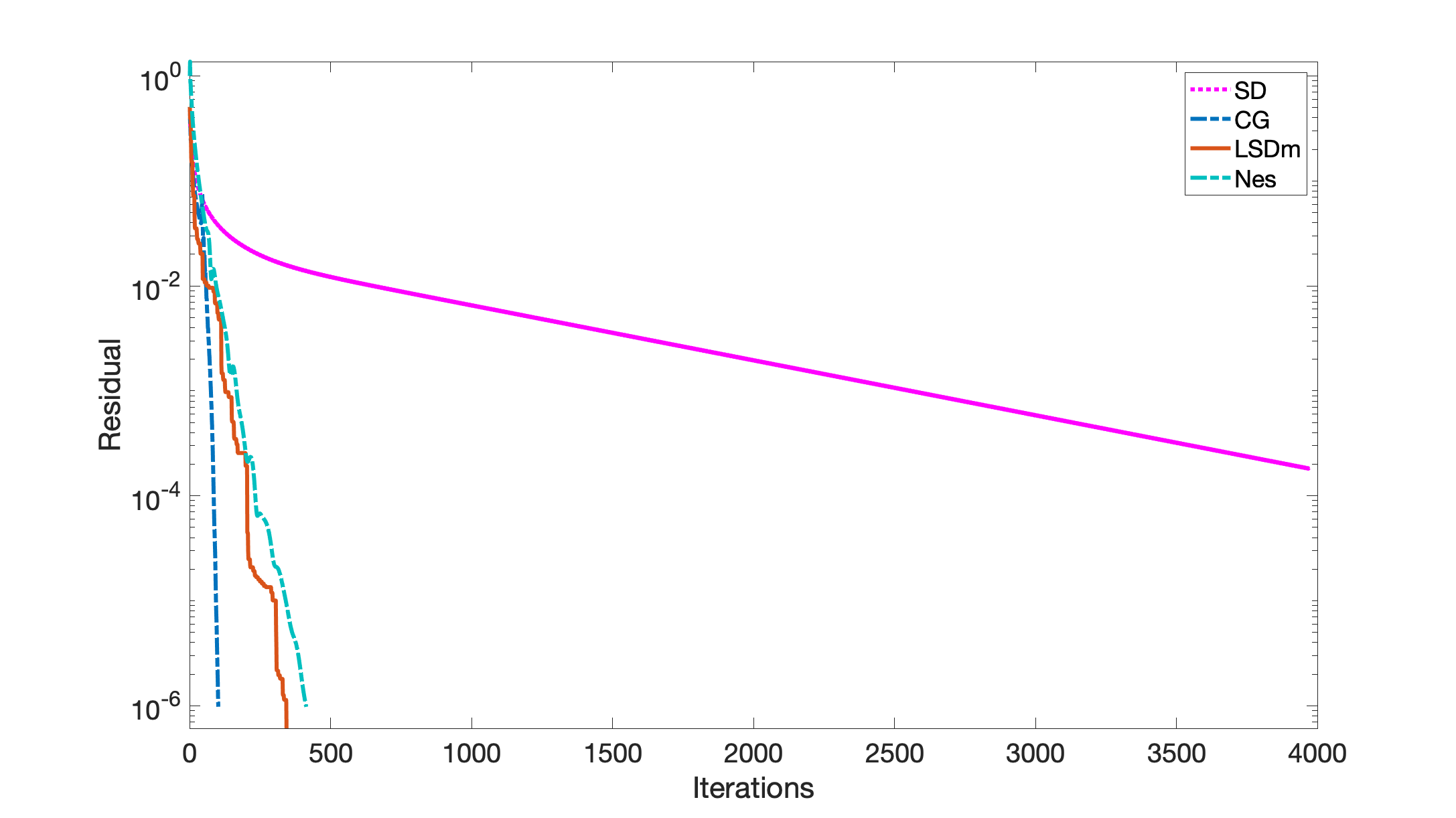}
\caption{Convergence behaviour of monotonized LSD (LSDm) as well as SD, conjugate gradient (CG)
and Nesterov's (Nes) for the model Poisson problem  with $n = 63^2$.
The errors $\| \rr_k \|$ 
are displayed as a function of iteration counter $k$. 
\label{fig:lsd_nesr}}
 \end{center}
\end{figure}

\begin{figure}[htb]
 \begin{center}
 \leavevmode
 \includegraphics[scale=0.17]{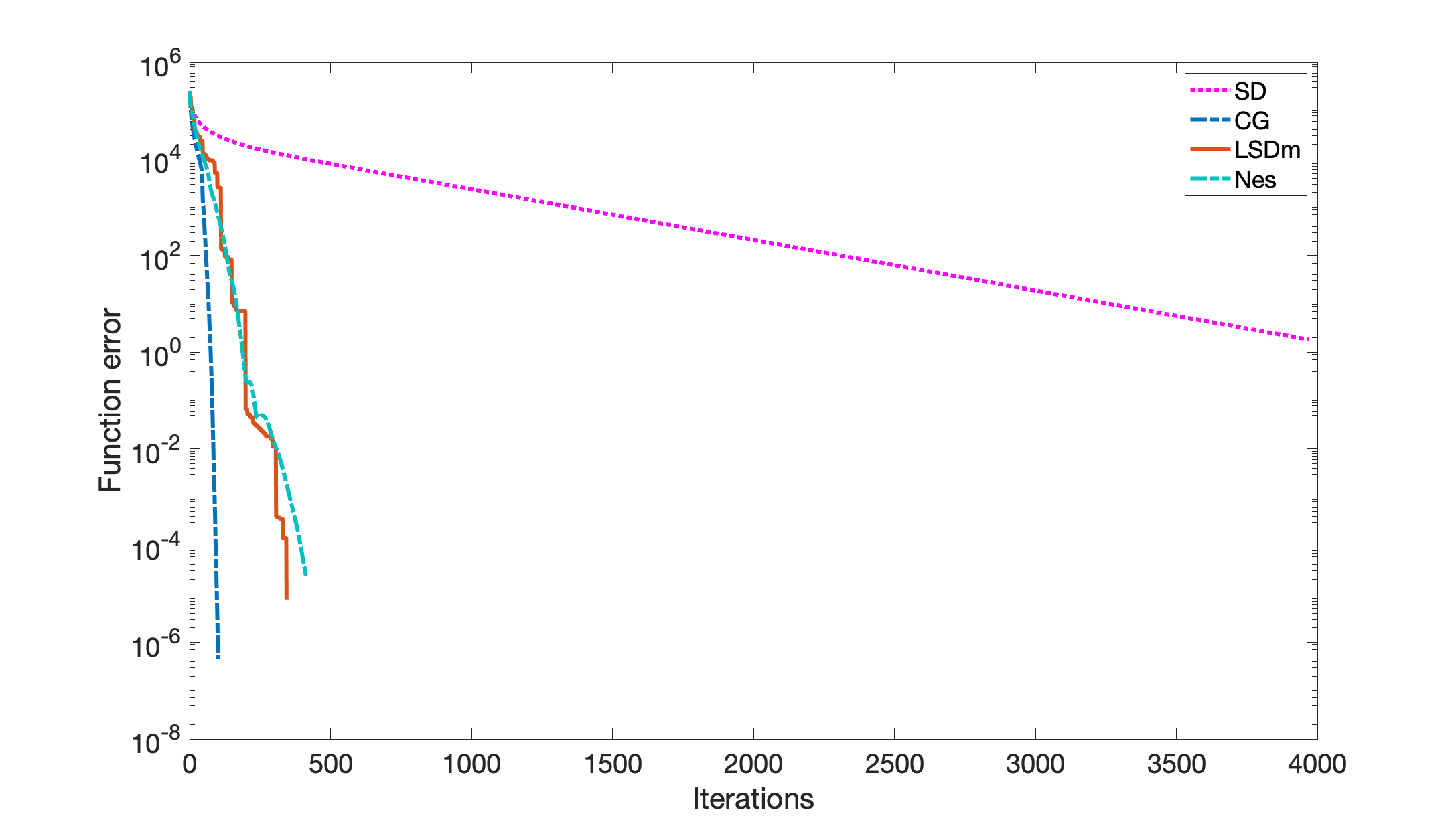}
\caption{Convergence behaviour of monotonized LSD (LSDm) as well as SD, conjugate gradient (CG)
and Nesterov's (Nes) for the model Poisson problem  with $n = 63^2$.
The errors $f( \xx_k) - f(\xx^*)$ 
are displayed as a function of iteration counter $k$. 
\label{fig:lsd_nesf}}
 \end{center}
\end{figure}

For the method~\eqref{9a5} there is the question of determining appropriate step sizes $\beta_k$ and $\alpha_k$.
If we keep both step sizes constant, then the iteration is stationary and standard techniques (involving the largest
and perhaps also smallest eigenvalues of the Hessian) apply.
For the calculations in Figures~\ref{fig:lsd_nesr} and~\ref{fig:lsd_nesf} we fixed the momentum parameter $\beta = 0.95$
\uri{(after some trial-and-error)} 
and calculated $\alpha_k$ by the SD formula~\eqref{5.6} with $\rr_k$ replaced by $\bb - A\yy_{k+1}$. 


In the present convex quadratic optimization context, all four methods depicted in these figures are Krylov-subspace methods  
generating iterates  in the same subspace at each $k$. So, the superiority of CG is a foregone conclusion,
as it is optimal on each such subspace (see, e.g.,~\cite{greenbaum}).
More interestingly, note the closeness of the LSD method and the accelerated gradient descent method~\eqref{9a5},
both being far better than SD and not much worse than CG. This suggests that the monotonized lagged steepest descent method
also satisfies the ``optimal'' estimate~\eqref{5.8}.\\

\paragraph{\bf Conclusion}
This longer case study, which involves some
new work, demonstrates a situation where the continuous 
differential equation limit~\eqref{11} misses out on some interesting and unusual action that occurs at the discrete process level.


\section{Gradient-based penalties in image processing}
\label{sec:paradigm}

In this case study we consider image and surface processing problems, where the given problem is entirely
discrete. 
Instances include image and surface denoising, deblurring
and inpainting, data consolidation, morphing, and much more. We can often view this as an instance of data
fitting: the given data $d$ relates to an image or object surface (e.g., it can be a point cloud),
the required operation such as denoising is described by a {\em forward operator} $\phi$,
and we seek another image or surface $u$ in \uri{2D or 3D, respectively,} such that $\phi(u)$ fits the data $d$ to a desired extent.     
Often, however,  a differential term is added in an attempt to regularize and control the desired process,
and the resulting differential problem is then subsequently discretized in order  to compute a solution.

\uri{Note that we are deliberately using the notation $u$ both for a discrete image/object and for this image/object interpolated in such
a way that its gradient in two or three space variables, denoted $\nabla u$, may be considered and differentiated. It is then
discretized (or restricted) back to the space where the data $d$ lives.
The step size $h$ is a parameter of such a discretization, say the distance between neighbouring image data locations.} 
Here we examine the utility of \uri{such a} continuous model, \uri{which the computer never ``sees''}. 

To be concrete, let us consider the problem
\uri{\begin{eqnarray}
\min_u f(u), \quad {\rm where~} f(u) = \| \phi (u) - d \|_2^2 + \beta \int_\Omega | \nabla u |^s , \label{41}
\end{eqnarray}}%
which we will call {\em the paradigm}, for short.
Here $\Omega$ is the domain of the image or the 3D object and $\beta > 0$ is \uri{a fixed, noise-dependent}
regularization weight as in Figure~\ref{fig1}.
Indeed,~\eqref{41} is in the typical form of Tikhonov regularization for an inverse problem~\cite{ta,ehn1,kaso}.
However, we stress that our concern here is with problems where $\phi$ is not a discretization of a continuous structure,
and the only source of differential terms, whether explicit or inverted, is the regularization term.
As such there is some cause for minor mathematical discomfort, because the leading (or, highest order) term of the 
PDE in the boundary value problem that arises in
the variational form  comes from the regularization rather than the data fitting term in~\eqref{41}.

There is a huge volume of literature on such methods; see, for instance,~\cite{sapiro,chsh,of,chpo2016}, to name but a few.
The parameter $s$ typically equals $2$, representing significant smoothing \uri{(diffusion)}, 
or $1$  for {\em total variation} (TV),
\uri{which typically yields piecewise smoothing and is also referred to as {\em anisotropic diffusion}}.
Other values such as $0 < s < 1$  are not our concern here.
The Euler-Lagrange equations for \eqref{41} give a boundary value PDE, and the problem actually solved is a discretization of such a PDE.
 The gradient in \eqref{41} can be replaced by a higher order differential operator, which is again not the focus here.
 If $s = 1$, then some \uri{further} regularization of the regularizer is required, and in~\cite{ahh} we proposed a method for controlling
 a Huber switch between $s=1$ and $s=2$ for this purpose. 
 
 There are lots of situations where the paradigm is the thing to consider! 
 \uri{Indeed, it introduces additional information which the discrete $\phi (u)$ as such does not have
 that relates the values of the $u$-unknowns to their neighbouring values. This in turn
 potentially enables a significantly simpler optimization process for $f(u)$.}
 
 However, already around the turn of the 21st century researchers
 found experimentally that the differential term in~\eqref{41} can introduce side effects that pollute the \uri{visual} quality of the obtained result.
 A ``social bifurcation'' \uri{has} formed, where researchers and practitioners in computer graphics and image processing
\uri{have}  turned more towards purely data driven methods, while mathematicians and numerical analysts \uri{have} continued to explore the paradigm and other such methods
 due to their more coherent structure that in turn enables development of better and more complete theory.
 
 Around the year 2005, armed with the code from~\cite{ahh}, we started exploring ways to introduce
 the paradigm for various image processing applications in computer graphics. Here the visual quality of the results is what matters,
 and no theory can compensate for lack in \uri{this} regard.
 We considered surface triangle mesh denoising in applications that involve intrinsic texture \uri{and
 sharp corners~\cite{huas,huas2}.} Together with computer graphics and vision experts we then considered  consolidation of unorganized point clouds 
 to allow quality surface reconstruction~\cite{hlzac}, point-resampling in the presence of edges~\cite{hwgcaz},
 and image tele-registration and structure-driven completion~\cite{hyglcac2013}.
 Unfortunately, while working on these projects our version of the paradigm was consistently beaten by other, discrete approaches,
 and there is no continuous model left in the final, \uri{published} version of any of
 \uri{\cite{huas,huas2,hlzac,hwgcaz,hyglcac2013}!} 
 
 Let us add a few more details regarding the research just mentioned.
 \begin{itemize}
 \item A noisy triangle mesh representing a surface has noisy
 triangle vertices, which are vectors in 3D. Several researchers have considered forms of 
 anisotropic Laplace-Beltramy PDE for denoising surface triangle meshes.
 But in \cite{huas} we learned (not before failing to extend the image denoising method of~\cite{tnv} to surfaces)
  that the most effective way to denoise such a vertex is to find the normal direction to the surface
 there and denoise in that specific direction. To find such a normal, the standard technique is to perform a principal component analysis (PCA)
 on a group of neighboring vertices. This process, described in detail in~\cite{huas}, is necessarily geometrically {\em local}, in marked
 difference to the {\em global} nature of an expression such as the penalty integral in~\eqref{41}.
 The obtained visual results are better, the complexity of the algorithm is linear in the number of unknowns(!), and      
 the corresponding {\sc Matlab} code has less than 50 lines. We feel that, unfortunately,
 one just can't practically do all that with anisotropic Laplace-Beltramy.
 \item
 In~\cite{hwgcaz} there is a construction of normals for a consolidated point cloud that represents a surface with an edge.
 Our algorithm is shown to perform better than  applying an $\ell_1$ regularization to the crude PCA normals
 (which is related to TV).
 \item
 The tele-registration algorithm in~\cite{hyglcac2013} involves a step to find salient curves in an image.
 This is a task for TV regularization, but the algorithm of~\cite{ffls2008} proved to be more robust.
 Further, towards the end of the image completion procedure there 
 the algorithm and software of~\cite{barnes2009,darabi2012}, which do not employ any \uri{direct} version of the paradigm,
 were used for  inpainting.
 \end{itemize}

 In~\cite{ashu18} we have summarized the pros and cons of the continuous paradigm~\eqref{41} vs discrete processes of
 data-driven reconstruction.
 
  
 \paragraph{\bf Conclusion} 
 The essential  advantage in using the paradigm is that one can often obtain a more solid theoretical backing to algorithms 
 \uri{that are occasionally more clearly derived based on solid principles} 
 from this global and generic point of view. The essential disadvantage, however, is that this approach leads to algorithms that 
 could be outperformed by more brute force techniques,  as described for the examples above. 
 This is so especially if (i) the forward operator $\phi$ is simple, and (ii) the data $d$ is of high \uri{visual} quality.


\section{Conclusions}
\label{sec:conclusions}

Recent years have seen a welcome development where knowledge and expertise  that have been gained over decades
in modelling and computationally solving differential equations
find use in the context of discrete problem areas such as optimization, \uri{game theory}, graphics and image processing.
The author, who has dealt with the numerical solution of differential equations most of his adult life
while being a member of a computer science department, finds special reasons to rejoice.
At the same time, we advocate balance and controlled euphoria. 
It is crucial to examine with an open mind the utility of such approaches and specific methods, 
from perspectives that include applicability,  relative advantage, and
extensions in the context of the task at hand.

We have demonstrated in a sequence of case studies that employing differential equations in this way can
be at times very useful, at times limiting, and at times just a matter of comfort.
One example of such ``comfort'', not mentioned hitherto, is the use of global properties of dynamical systems, such as
conservation laws, \uri{symplecticity} and reversibility, at the discrete level without requiring 
\uri{extreme} closeness of discrete and continuous model solutions.  
We advocate to become involved in this fascinating research avenue, as we ourselves are, while 
exercising practical caution.

\uri{We have deliberately avoided any significant discussion of stochastic optimization and stochastic differential equation methods
and models, as well as randomized algorithms. 
Our point is not to ignore these exciting areas, but rather to emphasize
that significant issues arise more basically even in their absence.}
   

\subsection*{Acknowledgment}
I thank Drs Jelena Diakonikolas, Arieh Iserles, Alexander Madureira, Fred Roosta, Lars Ruthotto, Mark Schmidt
and Jorge Zubelli for fruitful discussions.  
During several preparation stages of this paper I enjoyed the hospitality of IMPA, Rio de Janeiro. 



\bibliographystyle{plain}
\bibliography{biblio}

\begin{thebibliography}{10}

\bibitem{akaike}
H.~Akaike.
\newblock On a successive transformation of probability distribution and its
  application to the analysis of the optimum gradient method.
\newblock {\em Ann. Inst. Stat. Math. Tokyo}, 11:1--16, 1959.

\bibitem{a2}
U.~Ascher.
\newblock On symmetric schemes and differential-algebraic equations.
\newblock {\em SIAM J. Scient. Comput.}, 10:937--949, 1989.

\bibitem{colsys81}
U.~Ascher, J.~Christiansen, and R.D. Russell.
\newblock Collocation software for boundary-value {ODE}s.
\newblock {\em ACM Trans. Math. Software}, 7(2):209--222, 1981.

\bibitem{agbook}
U.~Ascher and C.~Greif.
\newblock {\em A First Course in Numerical Methods}.
\newblock SIAM, Philadelphia, PA, 2011.

\bibitem{ahh}
U.~Ascher, E.~Haber, and H.~Huang.
\newblock On effective methods for implicit piecewise smooth surface recovery.
\newblock {\em SIAM J. Sci. Comput.}, 28:339--358, 2006.

\bibitem{ashudo}
U.~Ascher, H.~Huang, and K.~van~den Doel.
\newblock Artificial time integration.
\newblock {\em BIT}, 47:3--25, 2007.

\bibitem{amr}
U.~Ascher, R.~Mattheij, and R.~Russell.
\newblock {\em Numerical Solution of Boundary Value Problems for Ordinary
  Differential Equations}.
\newblock SIAM, Philadelphia, 1995.

\bibitem{apbook}
U.~Ascher and L.~Petzold.
\newblock {\em Computer Methods for Ordinary Differential Equations and
  Differential-Algebraic Equations}.
\newblock SIAM, Philadelphia, PA, 1998.

\bibitem{ar2}
U.~Ascher and S.~Reich.
\newblock The midpoint scheme and variants for hamiltonian systems: advantages
  and pitfalls.
\newblock {\em SIAM J. Scient. Comput.}, 21:1045--1065, 1999.

\bibitem{ashu18}
U.M. Ascher and H.~Huang.
\newblock Numerical analysis in visual computing : What we can learn from each
  other.
\newblock {\em Vietnam J. Math}, pages DOI: 10.1007/s10013--018--0299--6, 2018.

\bibitem{bawi}
D.~Baraff and A.~Witkin.
\newblock Large steps in cloth simulation.
\newblock In {\em SIGGRAPH}, pages 43--54. ACM, 1998.

\bibitem{barnes2009}
C.~Barnes, E.~Shechtman, A.~Finkelstein, and D.~Goldman.
\newblock Patchmatch: a randomized correspondence algorithm for structural
  image editing.
\newblock {\em ACM Transactions on Graphics (SIGGRAPH)}, 27(3):67:1--67:10,
  2009.

\bibitem{babo}
J.~Barzilai and J.~Borwein.
\newblock Two point step size gradient methods.
\newblock {\em IMA J. Num. Anal.}, 8:141--148, 1988.

\bibitem{beck2017}
A.~Beck.
\newblock {\em First-Order Methods in Optimization}.
\newblock SIAM, 2017.

\bibitem{betancourt2018}
M.~Betancourt, M.~Jordan, and A.~Wilson.
\newblock On symplectic optimization.
\newblock {\em IEEE Trans. Visualization and Computer Graphics}, 2018.
\newblock arXiv1802.03653v2.

\bibitem{burger13}
M.~Burger, M.~Di Francesco, P.A. Markowich, and M.T. Wolfram.
\newblock On a mean field game optimal control approach modeling fast exit
  scenarios in human crowds.
\newblock In {\em IEEE Conference on Decision and Control}, volume~52, pages
  3128--3133, 2003.

\bibitem{care}
D.~Calvetti and L.~Reichel.
\newblock Lanczos-based exponential filtering for discrete ill-posed problems.
\newblock {\em Numer. Algorithms}, 29:45--65, 2002.

\bibitem{carezh}
D.~Calvetti, L.~Reichel, and Q.~Zhang.
\newblock Iterative exponential filtering for large discrete ill-posed
  problems.
\newblock {\em Numer. Math.}, 83:535--556, 1999.

\bibitem{chpo2016}
A.~Chambolle and T.~Pock.
\newblock An introduction to continuous optimization for imaging.
\newblock {\em Acta Numerica}, 25(161):161--319, 2016.
\newblock http://doi.org/10.1017/S096249291600009X.

\bibitem{chsh}
T.~Chan and J.~Shen.
\newblock {\em Image Processing and Analysis: Variational, PDE, Wavelet and
  Stochastic Methods}.
\newblock SIAM, 2005.

\bibitem{kale}
D.~Chen, D.~Levin, L.~Matusic, and D.~Kaufman.
\newblock Dynamics-aware numerical coarsening for fabrication design.
\newblock {\em ACM Trans. Graphics}, 36(4), 2017.

\bibitem{chorin}
A.J. Chorin.
\newblock Numerical solution of the navier-stokes equations.
\newblock {\em Math. Comp.}, 22:745--762, 1968.

\bibitem{chorin79}
A.J. Chorin and J.E. Marsden.
\newblock {\em A Mathematical Introduction to Fluid Mechanics}.
\newblock Springer, 1993.
\newblock 3rd ed.

\bibitem{darabi2012}
S.~Darabi, E.~Shechtman, C.~Barnes, D.~Goldman, and P.~Sen.
\newblock Image melding: Combining inconsistent images using patch-based
  synthesis.
\newblock {\em ACM Transactions on Graphics (SIGGRAPH)}, 31(4):82:1--82:10,
  2012.

\bibitem{deuflhard04}
P.~Deuflhard.
\newblock {\em Newton's Method for Nonlinear Problems}.
\newblock Springer, 2004.

\bibitem{dior19}
J.~Diakonikolas and L.~Orecchia.
\newblock The approximate duality gap technique: A unified theory of
  first-order methods.
\newblock {\em SIAM J. Optimiz.}, 29:660--689, 2019.

\bibitem{doas1}
K.~van~den Doel and U.~Ascher.
\newblock Dynamic level set regularization for large distributed parameter
  estimation problems.
\newblock {\em Inverse Problems}, 23:1271--1288, 2007.

\bibitem{doeasc}
K.~van~den Doel and U.~Ascher.
\newblock The chaotic nature of faster gradient descent methods.
\newblock {\em J. Scient. Comput.}, 51:560--581, 2011.

\bibitem{ehn1}
H.~W. Engl, M.~Hanke, and A.~Neubauer.
\newblock {\em Regularization of Inverse Problems}.
\newblock Kluwer, 1996.

\bibitem{ffls2008}
Z.~Farbman, R.~Fattal, D~Lischinski, and R.~Szeliski.
\newblock Edge-preserving decompositions for multi-scale tone and detail
  manipulation.
\newblock {\em ACM Transactions on Graphics (SIGGRAPH)}, 27(3):67:1--67:10,
  2008.

\bibitem{gomez14}
D.A. Gomez.
\newblock Mean field games models?a brief survey.
\newblock {\em Dynamic Games and Applications}, 4(2):110--154, 2014.

\bibitem{gomez13}
D.A. Gomez, J.~Mohr, and R.R. Souza.
\newblock Continuous time finite state mean field games.
\newblock {\em Applied Mathematics \& Optimization}, 68(1):99--143, 2013.

\bibitem{Goodfellow-et-al-2016}
I.~Goodfellow, Y.~Bengio, and A.~Courville.
\newblock {\em Deep Learning}.
\newblock MIT Press, 2016.
\newblock \url{http://www.deeplearningbook.org}.

\bibitem{greenbaum}
A.~Greenbaum.
\newblock {\em Iterative Methods for Solving Linear Systems}.
\newblock SIAM, Philadelphia, PA, 1997.

\bibitem{hlw}
E.~Hairer, C.~Lubich, and G.~Wanner.
\newblock {\em Geometric Numerical Integration}.
\newblock Springer, 2002.

\bibitem{hw}
E.~Hairer and G.~Wanner.
\newblock {\em Solving Ordinary Differential Equations II: Stiff and
  Differential-Algebraic Problems}.
\newblock Springer, 1996.
\newblock 2nd Edition.

\bibitem{huas}
H.~Huang and U.~Ascher.
\newblock Fast denoising of surface meshes with intrinsic texture.
\newblock {\em Inverse Problems}, 24 (3):034003, 2008.

\bibitem{huas2}
H.~Huang and U.~Ascher.
\newblock Surface mesh smoothing, regularization and feature detection.
\newblock {\em SIAM J. Scient. Comput.}, 31:74--93, 2008.

\bibitem{hlzac}
H.~Huang, D.~Li, R.~Zhang, U.~Ascher, and D.~Cohen-Or.
\newblock Consolidation of unorganized point clouds for surface reconstruction.
\newblock {\em ACM TOG (SIGGRAPH Asia)}, 29(5), 2009.

\bibitem{hwgcaz}
H.~Huang, S.~Wu, M.~Gong, D.~Cohen-Or, U.~Ascher, and H.~Zhang.
\newblock Edge-aware point set resampling.
\newblock {\em ACM TOG}, 32(1), 2013.

\bibitem{hyglcac2013}
H~Huang, K~Yin, M~Gong, D~Lischinski, D~Cohen-Or, UM~Ascher, and B~Chen.
\newblock ``mind the gap": tele-registration for structure-driven image
  completion.
\newblock {\em ACM Transactions on Graphics}, 32(6):174:1--174:10, 2013.

\bibitem{kaso}
J.~Kaipo and E.~Somersalo.
\newblock {\em Statistical and Computational Inverse Problems}.
\newblock Springer, 2005.

\bibitem{kreiss72}
H.-O. Kreiss.
\newblock Centered difference approximation to singular systems of odes.
\newblock {\em Symposia Mathematica X}, 1972.
\newblock Inst. Nazionalle di Alta Math.

\bibitem{ll07}
J.~Lasry and P.~Lions.
\newblock Mean field games.
\newblock {\em Jpn. J. Math.}, 2(1):229--260, 2007.

\bibitem{hmc06}
R.~Malhame M.~Huang and P.~Caines.
\newblock Large population stochastic dynamic games: Closed-loop mckean?vlasov
  systems and the nash certainty equivalence principle.
\newblock {\em Comm. Info. Sys.}, 6(3):221--252, 2006.

\bibitem{mrs90}
M.~Markowich, C.~Ringhofer, and C.~Schmeiser.
\newblock {\em Semi-conductor equations}.
\newblock Springer, 1990.

\bibitem{nesterov1983method}
Y.~Nesterov.
\newblock A method of solving a convex programming problem with convergence
  rate o(1k2).
\newblock In {\em Doklady Akademii Nauk}, pages 543--547, 1983.

\bibitem{nesterov18}
Y.~Nesterov.
\newblock {\em Lectures on convex optimization}.
\newblock Springer, 2018.

\bibitem{nw}
J.~Nocedal and S.~Wright.
\newblock {\em Numerical Optimization}.
\newblock New York: Springer, 1999.

\bibitem{of}
S.~Osher and R.~Fedkiw.
\newblock {\em Level Set Methods and Dynamic Implicit Surfaces}.
\newblock Springer, 2003.

\bibitem{polyak64}
B.~T. Polyak.
\newblock Some methods of speeding up the convergence of iteration methods.
\newblock {\em USSR Comput. Math. \& Math. Phys.}, 4:1--17, 1964.

\bibitem{rasv}
M.~Raydan and B.~Svaiter.
\newblock Relaxed steepest descent and {C}auchy-{B}arzilai-{B}orwein method.
\newblock {\em Comp. Optimization Applic.}, 21:155--167, 2002.

\bibitem{rodoas}
F.~Roosta-Khorasani, K.~van~den Doel, and U.~Ascher.
\newblock Stochastic algorithms for inverse problems involving {PDE}s and many
  measurements.
\newblock {\em SIAM J. Scient. Comput.}, 2014.

\bibitem{ruha19}
L.~Ruthotto and E.~Haber.
\newblock Deep neural networks motivated by partial differential equations.
\newblock {\em J. Mathematical Imaging and Vision}, 2019.

\bibitem{sapiro}
G.~Sapiro.
\newblock {\em Geometric Partial Differential Equations and Image Analysis}.
\newblock Cambridge, 2001.

\bibitem{su2014}
W.~Su, S.~Boyd, and E.~Candes.
\newblock A differential equation for modelling nesterov's accelerated gradient
  method.
\newblock {\em Advances in Neural Information Processing Systems (NIPS)}, 27,
  2014.

\bibitem{su2016}
W.~Su, S.~Boyd, and E.~Candes.
\newblock A differential equation for modelling nesterov's accelerated gradient
  method: Theory and insights.
\newblock {\em J. Machine Learning Research}, 17(153):1--43, 2016.

\bibitem{tnv}
E.~Tadmor, S.~Nezzar, and L.~Vese.
\newblock A multiscale image representation using hierarchical ({BV},${L}^2$)
  decompositions.
\newblock {\em SIAM J. Multiscale Model. Simul.}, 2:554--579, 2004.

\bibitem{ta}
A.~N. Tikhonov and V.~Ya. Arsenin.
\newblock {\em Methods for Solving Ill-posed Problems}.
\newblock John Wiley and Sons, Inc., 1977.

\bibitem{wanner10}
G.~Wanner.
\newblock {K}epler, {N}ewton and numerical analysis.
\newblock {\em Acta Numerica}, 19:561--598, 2010.

\bibitem{wibisono2016}
A.~Wibisono, A.~Wilson, and M.~Jordan.
\newblock A variational perspective on accelerated methods in optimization.
\newblock {\em Proc. Nat. Academy Science}, 113(47):E7351--E7358, 2016.

\bibitem{zangwill69}
Willard~I. Zangwill.
\newblock {\em Nonlinear programming: a unified approach}.
\newblock Prentice-Hall, Inc., Englewood Cliffs, N.J., 1969.

\bibitem{zhang18}
J.~Zhang, A.~Mokhtari, S.~Sra, and A.~Jadbabai.
\newblock Direct runge-kutta discretization achieves acceleration.
\newblock {\em NeurIPS}, 2018.
\newblock arXiv:1805.00521.

\end{thebibliography}

\end{document}